\documentclass[11pt,reqno,a4paper,twoside]{article}

\usepackage{amsmath,amsthm,amstext,amscd,amssymb,euscript,mathrsfs, url}
\usepackage{calrsfs}
\usepackage{epsfig}

\newcommand{\eec}{\color{black}}
\newcommand{\bwc}{\color{black}}
\newcommand{\ewc}{\color{black}}

\newcommand{\Z}{\mathbb Z}
\newcommand{\R}{\mathbb R}
\newcommand{\N}{\mathbb N}

\newcommand{\E}{\mathbb E}
\newcommand{\Zd}{\mathbb Z^d}

\renewcommand{\phi}{\varphi}
\newcommand{\epsi}{\ensuremath{\epsilon}}

\newcommand{\si}{\ensuremath{\sigma}}

\newcommand{\pee}{\ensuremath{\mathbb{P}}}

\def\1{{\mathchoice {\rm 1\mskip-4mu l} {\rm 1\mskip-4mu l}
{\rm 1\mskip-4.5mu l} {\rm 1\mskip-5mu l}}}

\newtheorem{theorem}{{\small T}{\scriptsize HEOREM}}[section]
\newtheorem{corollary}{{\bf{\small C}{\scriptsize OROLLARY}}}[section]
\newtheorem{proposition}{{\bf{\small P}{\scriptsize ROPOSITION}}}[section]
\newtheorem{lemma}{{\bf{\small L}{\scriptsize EMMA}}}[section]
\newtheorem{remark}{{\bf{\small R}{\scriptsize EMARK}}}[section]
\newtheorem{definition}{{\bf{\small D}{\scriptsize EFINITION}}}[section]

\renewenvironment{proof}[1]
{\noindent{{\bf{\small{ P}{\scriptsize ROOF}}}.}\hspace{0.1cm} #1} {$\;\qed$\newline}

\newcommand{\beq}{\begin{eqnarray}}
\newcommand{\eeq}{\end{eqnarray}}

\newcommand{\ba}{\begin{align*}}
\newcommand{\ea}{\end{align*}}

\newcommand{\be}{\begin{equation}}
\newcommand{\ee}{\end{equation}}

\newcommand{\bl}{\begin{lemma}}
\newcommand{\el}{\end{lemma}}

\newcommand{\br}{\begin{remark}}
\newcommand{\er}{\end{remark}}

\newcommand{\bt}{\begin{theorem}}
\newcommand{\et}{\end{theorem}}

\newcommand{\bd}{\begin{definition}}
\newcommand{\ed}{\end{definition}}

\newcommand{\bp}{\begin{proposition}}
\newcommand{\ep}{\end{proposition}}

\newcommand{\bc}{\begin{corollary}}
\newcommand{\ec}{\end{corollary}}

\newcommand{\bpr}{\begin{proof}}
\newcommand{\epr}{\end{proof}}

\newcommand{\bi}{\begin{itemize}}
\newcommand{\ei}{\end{itemize}}

\newcommand{\ben}{\begin{enumerate}}
\newcommand{\een}{\end{enumerate}}

\newcommand{\caR}{{\mathcal R}}
\newcommand{\caS}{{\mathcal S}}
\newcommand{\caT}{{\mathcal T}}

\usepackage{graphics,color,calc,graphicx}

\usepackage[a4paper,vmargin={3.5cm,3.5cm},hmargin={2.5cm,2.5cm},bottom=30mm,
marginparwidth=50pt,marginparsep=1mm]{geometry}

\def\1{{\mathchoice {\rm 1\mskip-4mu l} {\rm 1\mskip-4mu l}
{\rm 1\mskip-4.5mu l} {\rm 1\mskip-5mu l}}}

\begin{document}
\title{Non-criticality criteria for Abelian sandpile models with sources and sinks}
\author{
Frank Redig$^{\textup{{\tiny(a)}}}$, Wioletta M. Ruszel$^{\textup{{\tiny(a)}}}$ and
Ellen Saada$^{\textup{{\tiny(b)}}}$
\\
{\small $^{\textup{(a)}}$
Delft Institute of Applied Mathematics,}\\
{\small Technische Universiteit Delft}\\
{\small Van Mourik  Broekmanweg 6, 2628 XE Delft, Nederland}
\\
{\small $^{\textup{(b)}}$}
\\
{\small CNRS, UMR 8145, Laboratoire MAP5,}\\
{\small
Universit\'{e} Paris Descartes, Sorbonne Paris Cit\'{e},} \\
{\small 45, rue des Saints P\`{e}res 75270 Paris Cedex 06, France}\\
}
\maketitle
\begin{abstract}
We prove that the Abelian sandpile model on a random binary and binomial tree, as
introduced in [F. Redig, W. M. Ruszel, E. Saada, \textit{The Abelian sandpile model on a random tree},
\textit{J. Stat. Phys.} {\bf 147},  653--677, (2012)], is not critical for all branching probabilities $p<1$;
by estimating the tail of the annealed survival time of a random
walk on the binary tree
with randomly placed traps, we  obtain  some more information
about the exponential tail of the avalanche radius.
Next we study the sandpile model on $\Zd$ with some additional dissipative sites: 
we provide examples and sufficient conditions for non-criticality;
we also make a connection with the parabolic Anderson model.
Finally we initiate the study of the sandpile model with both 
sources and sinks and give  a sufficient condition for 
non-criticality in the presence of a finite number of sources, 
using a connection with the homogeneous pinning model.
\end{abstract}
\section{Introduction}
The Abelian sandpile model was introduced by \cite{btw} and \cite{dharb}
as a toy model displaying self-organized criticality. 
This means that in the thermodynamic limit, the sandpile model
has features of models of statistical physics at the critical point,
such as power law decay of avalanches, but without any fine-tuning of
parameters.
In short, an Abelian sandpile model is a discrete dynamical 
system defined as follows. Assign to each vertex of some 
finite graph a discrete variable (number of particles) 
representing a height, and call some special vertices sinks.
 At each time, we add an extra particle to the system uniformly 
 at random. If the resulting height exceeds a certain threshold, then the 
 vertex topples by distributing the particles among its neighbours. 
 This toppling can cause other vertices to topple which can lead 
 to an avalanche. The self-organized criticality comes into play 
 in the fact that the avalanche distribution obeys in certain cases a power-law. 
In order to stabilize, the sandpile model needs sinks, which in the standard setting
are associated with the boundary. On the contrary, if grains are lost upon
toppling the bulk, we call the model dissipative. It is well-known that this leads to
non-criticality (i.e., exponential decay of avalanche sizes). In this paper, we are interested
in how much dissipation is needed in order to lose criticality.
In the language of random walks, this amounts to create traps (associated to dissipative sites) and the expected
avalanche size corresponds to an expected survival time, which depends on the configuration of traps.
Second, we also introduce  sites where upon toppling mass is gained (source sites).
In the random walk picture, this corresponds to branching and the combination of sources and sinks
to random walk with branching and killing. We associate non-criticality with finite expected avalanche size which then, in the setting
of models with sources and sinks, corresponds to the finite expected total mass, when the random walk with killing and branching
is started with unit mass.
   
In \cite{rrs}, we introduced the Abelian sandpile model on a 
random binary (or binomial) tree and
proved that in a small region of the supercritical regime of the  associated 
branching process, the model is not critical: that is,
avalanche sizes decay exponentially. 
For the full Bethe lattice, avalanche sizes decay as a power law (see \cite{dhar}), 
and so it remained an open question whether the absence of criticality persists 
in the whole supercritical
regime, except for the degenerate case of the full Bethe lattice. 
In the present paper, we study this question, thereby solving  an  open issue of \cite{rrs}. 

In relation to this problem, we also tackle the more general question: how much dissipation 
can be added to the sandpile model in order
to make it non-critical?
 Therefore  we consider  various other settings and examples
  where we obtain criteria for non-criticality, such as $\Zd$ with
 a (possible random) set of  dissipative sites. 
Furthermore  we also start the study of the sandpile model with both sources and sinks.
Assuming that the height variables  do not grow indefinitely, i.e., 
if the system is stabilizable, we show that a finite number of source sites 
surrounded by enough dissipative sites produces a non-critical system. 
Even this at first sight intuitively obvious fact turns out to be non-trivial to prove; 
indeed, in the presence of source sites, large deviations come into play, and 
therefore, the simple intuition that the system is non-critical when there are 
more dissipative sites than source sites is wrong. For example, the system with one source 
site and all other sites ``standard'' is not stable.

The first result in this direction obtained in our paper also creates 
a link between the sandpile and the parabolic Anderson model,  which is a model of a random walk in a random potential landscape where sites with negative potential are interpreted as killing and positive potential as branching, 
(see \cite{konig} for an up-to-date introduction to this model): 
  we show that  the non-criticality of the sandpile  model with randomly placed 
 sinks and sources is equivalent to the finiteness of the first moment 
 of the total mass in the corresponding parabolic Anderson model. The source 
 sites act as sites where the random walk (of the corresponding parabolic 
 Anderson model) is branching, whereas the dissipative sites act as 
 (possibly soft) traps. Therefore, whether or not the expected total mass 
 has a first moment is related to random walk local time large deviation properties. \\

The rest of the paper is organized as follows. In section \ref{sec:survive} 
we prove the exponential tail of the \bwc probability of survival \ewc of a random walk on a $q$-ary tree with traps. 
In section \ref{sec:ranbin} we apply this result to prove non-criticality of 
the sandpile model on a random binomial tree and more general supercritical 
branching processes with bounded offspring distribution. In section \ref{sec:diss}
 we study the sandpile model on $\Zd$ and prove, guided by examples, sufficient 
 conditions for non-criticality in the presence of dissipative sites. In section 
 \ref{sec:disssource} we consider the model with sources and sinks and prove a 
 first result of non-criticality in that context, making use of a  connection 
 with the homogeneous pinning model. 
%
%
\section{\bwc Probability of survival \ewc of a random walk on the $q$-ary tree with traps}\label{sec:survive}
In this section we consider a random walk on a $q$-ary tree 
with randomly placed traps and prove that \bwc the probability of \ewc its survival decays exponentially.
There is a large literature on the \bwc probability of \ewc survival  asymptotics for the walk on $\Zd$, 
which  is related  to the large deviations of the Wiener sausage
 (see e.g.  \cite{deuschelstroock} and \cite{dv}).
Surprisingly, we could not find the corresponding result for random walk on trees, 
although some results on trapped random walks on trees exist  (see e.g.  \cite{konig}). 
Intuitively,  the probability of survival should decay exponentially 
because the random walk is strongly transient;  hence, it has a range which 
typically grows linearly in the number of steps.
This intuition is exactly what we make rigorous in this section.

Consider a rooted infinite $q$-ary tree  $\caT_q$ with root  $o$; 
i.e. each vertex  $x \in \caT_q$ different from the root 
has degree $q+1$, and the root has degree $q$.  

For two vertices $x,y \in \caT_q$, we denote by $d(x,y)$ the graph distance 
between $x$ and $y$  on $\caT_q$. 
We place at every site $x\in \caT_q$
a trap with probability $p_x>0$, independently for different sites.
We call $\omega: \caT_q \to \{0,1\}$ a trap configuration. 
More precisely, we set $\omega(x)=1$ if $x$ is a trap with probability $p_x$, 
respectively $\omega(x)=0$ otherwise. Note that we allow the trapping probability 
to depend on the location $x$.

Denote by $(S_n)_{n\geq 0}$ a simple random walk on $\caT_q$.
This random walk  is killed upon hitting a trap with probability 1, and we call
$T(\omega)$ its survival time,
i.e.,
\be\label{traptime}
T(\omega)=\inf\{ n\in \N: \omega(S_n)=1\}.
\ee
 We denote by $\bf{P}_{o}(\cdot)$ \bwc the annealed law, \ewc the joint probability of  the random walk  
 $(S_n)_{n\geq 0}$ on 
$\mathcal{T}_q$ starting at $o$ together with the trap configuration $\omega$ 
conditioned on the root not being a trap, and we denote by 
$\bf{E}_o(\cdot)$ the corresponding expectation. We condition on the fact 
that the root is not a trap in order to make the link  to branching processes 
in  Section \ref{sec:ranbin}.
In Proposition \ref{expobound} below, we prove that the probability of survival 
decays exponentially in the annealed
setting.
This is in contrast to the analogous problem
on $\Zd$, where the optimal strategy to survive  for a long time is to
stay for a long time in a trapless ball around the origin \cite{konig}.
\bp\label{expobound}
\bwc Assume that  there exists $p>0$ such that the trapping probabilities 
$(p_x)_{x\in \mathcal{T}_q}$ satisfy  $p_x>p$  for all $x$  \ewc.
Then there exists  $c=c(p,q)>0$ such that for all $n$ 
\be
\textbf{P} _o (T(\omega)>n)\leq e^{-cn}, 
\ee
and as a consequence all moments of $T$ exist.  
\ep
\bpr
We have  $S_0= o$; call $d(o,S_n)= :X_n$.
Then,  $(X_n)_{n\geq 0}$  is a one-dimensional random walk 
making a $+1$ step with probability $q/(q+1)$
and a $-1$ step with probability $1/(q+1)$, reflected at the origin.

Let us denote by $R_n(\omega)$ the range of the random walk $(S_n)_{n\geq 0}$, i.e.,
the number of points visited before time $n$.
Then we have that on the event $\{T(\omega)>n\}$, 
all the points  counted in $R_n(\omega)$ should be trapless so that
\begin{eqnarray}\label{bo}
\textbf{P}_o(T(\omega)>n)&\leq& \textbf{E}_o ((1-p)^{R_n(\omega)})
\nonumber\\
&\leq & \textbf{E}_o \left ((1-p)^{\epsilon n} \1_{\{R_n(\omega) 
\geq \epsilon n \}}\right ) + \textbf{E}_o( \1_{ \{R_n(\omega) <\epsilon n \}}) \ewc
\nonumber\\
&\leq &
(1-p)^{\epsi n} + \textbf{P}_o(R_n(\omega)<\epsilon n).
\end{eqnarray}
To estimate $\textbf{P}_o(R_n(\omega)<\epsilon n)$, realize that 
if $X_n \geq  \epsi n$ then necessarily
$R_n(\omega) \geq  \epsi n$. Hence if $R_n(\omega) < \epsilon n$, then 
also $X_n < \epsilon n$. Therefore, because the random walk  $(X_n)_{n\geq 0}$  has a drift
\[\frac{q}{q+1}- \frac{1}{q+1}= \frac{q-1}{q+1},\]
 by  applying   the Chernoff bound for $t>0$, we have
\[
\textbf{P}_o (X_n < \epsi n) \leq e^{\epsi n t} \left ( \frac{q}{1+q}e^{-t}
 + \frac{1}{q+1}e^{t}\right )^n.
\]
Choose 
\[t:=\frac{1}{2}\log \left ( \frac{2q}{1+\epsi} -q\right)>0\]
 then
there exists $\epsi:=\epsi(q)$ such that for all $c < \epsi$ and all $n$ we have that
\be\label{baba}
\textbf{P}_o(R_n (\omega)<\epsilon n)\leq \textbf{P}_o (X_n < \epsi n)  \leq  e^{-c n}
\ee
 \bigskip
which,  combined with \eqref{bo},  yields the claim.
\epr

Let us  consider   $\{S_n, 0\leq n\leq T(\omega)\}$,  i.e. 
 the random walk killed upon hitting a trap.
First note that for every given environment $\omega$ of traps, 
we can write the survival time $T(\omega)$ of the random walk as
\begin{equation}\label{def:eq_T}
T(\omega)= \sum_{ x \in \caT_q } \sum_{n=0}^{T(\omega)}  \1_{ \{S_n=x \}}.
\end{equation}
We denote by $\E^{\omega}_{o}(\cdot)$ the  quenched  expectation over the random walk 
 started at $o$ with a \bwc \textit{fixed} \ewc trap configuration $\omega$, and by $E(\cdot)$ the average over $\omega$. 
We have by taking the expectation first over the random walk and second over the traps,
\be\label{sumgreen}
 \mathbf{E}_o (T(\omega)) =   \mathbf{E}_o \left ( \sum_{x\in \caT_q} \sum_{n=0}^{T(\omega)} \1_{\{S_n=x\}}  \right ) = E \left (\sum_{x \in \caT_q} \mathcal{G}_{T(\omega)} (o, x) \right ),
\ee
where
\begin{equation}\label{def:eq_G}
\mathcal{G}_{T(\omega)} (o,x) =  \E^{\omega}_{o} \left ( \sum_{n=0}^{T(\omega)} \1_{ \{S_n=x\}} \right ) 
\end{equation}
is the Green's function of  the random walk $(S_n)_{n\geq 0}$  started at $o$, for a given trap configuration 
$\omega$. 

 The following annealed bound  will be useful in Section \ref{sec:ranbin}
to estimate avalanche diameter of the sandpile model. 
\bp\label{diam}
There exists a constant  $C:= C(p,q)>0$  such that
\be\label{diamgreen}
 E\left (\sum_{x: d(0,x)>n} \mathcal{G}_{T(\omega)} (o,x) \right ) \leq e^{-C n}.
\ee
\ep
\bpr
\begin{eqnarray}
\sum_{x: d(0,x)>n} \mathcal{G}_{T(\omega)} (o,x) & = &  
\sum_{x: d(0,x)>n} \E^{\omega}_{o} \left ( \sum_{n=0}^{T(\omega)} \1_{ \{S_n=x\}} \right ) \\
&=&
\sum_{x: d(0,x)>n} \mathbb{E}^{\omega}_o \left (\sum_{k=1}^{T(\omega)} \1_{\{S_k=x \}} 
\1_{\{T(\omega)>n\}}\right )
\nonumber\\
&=&
\mathbb{E}^{\omega}_o \left (\1_{\{T(\omega)>n\}} \sum_{x: d(0,x)>n}\sum_{k=1}^{T(\omega)} 
\1_{ \{S_k=x\}} \right )
\nonumber\\
& \leq &
 \mathbb{E}^{\omega}_o \left (\1_{\{T(\omega)>n\}} T (\omega) \right ),
\end{eqnarray}
where the first equality comes from \eqref{def:eq_G} and the inequality from  \eqref{def:eq_T}. 
As a consequence, taking the expectation over all the trap configurations provides
\[
\begin{split}
E \left(\sum_{x: d(0,x)>n} \mathcal{G}_{T(\omega)} (o, x)\right) & \leq E \left( \mathbb{E}^{\omega}_o (T(\omega) \1_{\{T(\omega)>n\}} )\right). 
\end{split}
\]
By  using first Cauchy-Schwartz inequality then Proposition \ref{expobound},  we have for some $C>0$,
\[
\textbf{E}_o \left(T(\omega) \1_{\{T(\omega)>n\}} \right) \leq \sqrt{\textbf{E}_o(T^2(\omega))}\sqrt{\textbf{P}_o(T(\omega) >n)} \leq e^{-C n}
\]
and the claim follows. 
\epr

\section{The sandpile model on the random binomial tree}\label{sec:ranbin}

In this section, we will first describe how we can construct a realization 
of a Galton-Watson branching process $\mathbb{T}^q$ starting from a single 
individual and having at most  $q$ descendants from  a realization of a  $q$-ary tree 
$\mathcal{T}_q$ with  some appropriate deletion, and second we will define 
the sandpile model on this realization.

We identify any vertex $x$ (apart from the root) of a $q$-ary tree 
$\mathcal{T}_q$   with a trap 
with probability $p>0$ and no trap with probability $1-p$ 
(thus there is no dependency on the location  $x\not= o$, but the root $o$ is trapless). 
This defines a trap configuration $\omega$. 
 If a vertex is a 
trap, then we delete it, as well as all its descendants and the edge to its parent. 
Note that the remaining tree is in distribution equal to a realization of a Galton-Watson 
branching process with offspring distribution given by $Bin(q,1-p)$. 
We will refer to this tree  as a random binomial tree with at most  
$q$ descendants. Remark that for this construction to hold we had to assume 
that the root is not a trap.

We now consider a sandpile model on a fixed realization of the binomial tree. 
We first define the model in finite volume and then give  its 
infinite volume construction. Standard references are  \cite{jarai, jarai1, redig} and \cite{redigo}. 

For some $n \in \mathbb{N}$,  denote by 
$\mathbb{T}^q_n=\{ x\in \mathbb{T}^q: d(o,x)\leq n \}$ the random binomial tree 
with root $o$ up to generation $n$. We assign to each 
vertex $x$  a height
$\eta_x \in \mathbb{N}$.
A  height configuration $\eta=(\eta_x)_{x\in \mathbb{T}_n^q}$ is {\it stable} 
if for all sites $x\in\mathbb{T}^q_n$, $\eta_x\in \{0,1,2,...,q\}$;  otherwise we call it {\it unstable}. 
If a site $x$ is unstable (i.e., if $\eta_x>q$), it topples:  this means that
it ejects  $q+1$ particles out of which it redistributes one particle to each of its neighbours. 
If there exists particles ejected that are not redistributed, then the site is said \textit{dissipative.}
Note that the root is always losing 1 particle upon toppling.
The effect of  toppling of  site $x$  on  the height configuration $\eta$ 
at some site $y \in \mathbb{T}^q_n$ is given by 
\[
(T_x\eta)_y = \eta_y- (\Delta^{\mathbb{T}_n^q})_{x,y}.
\]
 Here $\Delta^{\mathbb{T}_n^q}$, the
toppling matrix, is
the discrete Laplacian with Dirichlet boundary conditions, i.e.,
\be\label{top}
(\Delta^{\mathbb{T}_n^q})_{x,y}=
\begin{cases}
q+1\ & \text{if}\ x=y\in \mathbb{T}_n^q
\\
-1\ & \text{if}\ x\sim y \in \mathbb{T}_n^q
\\
0\ & \text{otherwise},
\end{cases}
\ee
where $x\sim y$ means $d(x,y)=1$.
 
We denote by  $\Omega_{n}$ the set of stable  configurations
on $\mathbb{T}_n^q$. We call $y$ a boundary site of $\mathbb{T}_n^q$ if $d(o,y)=n$. 
Note that in this model every site which has either non-maximal degree or is  a boundary site
 is dissipative.  
As a consequence, because not all vertices have maximal degree
with positive probability, there is a positive density
of dissipative sites outside the boundary,  which intuitively should lead to a non-critical model.
Note that this set-up, similar to what was done in \cite{rrs}, 
is in contrast with the set-up of \cite{jrs1}, where the random toppling matrix 
depends on the realization of the tree in such a way that
the only dissipative sites are the boundary sites, so that the model there is critical. \\

A toppling at $x\in\mathbb{T}_n^q$ in configuration $\eta$  is called \textit{legal} 
if $\eta_x > q$. A sequence of legal topplings is a composition $T_{x_\ell}\circ\ldots\circ T_{x_1} (\eta)$ such that
for all $k=1,\ldots,\ell$,  we have $x_k\in\mathbb{T}_n^q$ and the toppling at $x_k$ is legal
in $T_{x_{k-1}}\circ\ldots\circ T_{x_1} (\eta)$.  

We denote by $\caS_n$ the stabilization with toppling matrix \eqref{top}, and by
$a_{x,n}$ the addition operator defined on $\Omega_n$ via
\be
a_{x,n} (\eta)=\caS_n (\eta+\delta_x),
\ee
  where $\delta_x$ denotes a unit mass at $x$ and zero mass everywhere else. In other words, 
$\caS_n (\eta+\delta_x)$ is the unique stable configuration 
that arises from $\eta+\delta_x$ by a sequence of legal topplings.
Notice that both $\caS_n$ and $a_{x,n}$ are quenched, i.e., for a given realization
$\mathbb{T}^q$ of the random tree. Thus 
the model is Abelian so that the stabilization does not depend on the 
toppling order and indeed leads to a unique stable configuration.
The set of recurrent configurations corresponding to the addition operators
$a_{x,n}$
is denoted by $\caR_n$. The unique measure $\mu_n$ on $\Omega_n$ invariant under all addition operators $a_{x,n}$ 
is the uniform measure on the set of recurrent configurations
\[
\mu_n = \frac{1}{|\caR_n|}\sum_{\eta \in \caR_n} \delta_{\eta}.
\]

By Theorem 3  of \cite{jw12} an infinite volume 
sandpile measure $\mu=\lim_{n\rightarrow \infty} \mu_n$ on $\mathbb{T}^q$ exists because each tree 
in the wired uniform spanning forest on $\mathbb{T}^q$ has one 
end almost surely. See \cite{LP16} for the background on wired spanning forests.
The one end property for Galton-Watson trees with bounded degree 
distribution was proven in Theorem 7.2 in \cite{AL07}.  

For $\mu$-a.e. $\eta$, the addition operator
\[
a_x (\eta)=\lim_{n\to\infty} a_{x,n} (\eta)
\]
is well defined (see \cite{jarai}), where with a small abuse of the notation we denoted by $a_{x,n} (\eta)$ the concatenation of the addition operator applied in finite volume $\mathbb{T}^q_n$ and the identity, i.e. 
fixing $\eta$ outside $\mathbb{T}^q_n$, that is 
$a_{x,n} (\eta)= a_{x,n} (\eta|_{|\mathbb{T}^q_n}) \eta_{|(\mathbb{T}^q_n)^c}$. 
 
Moreover, we have
\be\label{bori}
a_x (\eta)= \eta+\delta_x- \sum_{y\in\mathbb{T}^q} (\Delta^{\mathbb{T}^q})_{x,y} N(x, y,\eta),
\ee
where $N(x,y,\eta)$ denotes the number of topplings needed at $y$ 
to stabilize $\eta+\delta_x$, also known as the odometer function, 
and $\Delta^{\mathbb{T}^q}$ is defined analogously to \eqref{top}, 
with $\mathbb{T}_n^q$ replaced by $\mathbb{T}^q$. 
\eec
We then define the avalanche at $x$ in configuration $\eta$ by
\be
Av (x,\eta)=\{ y: N(x, y,\eta)>0\}.
\ee
Integrating \eqref{bori} with respect to $\mu$ (we denote by $\E_{\mu}$ the corresponding expectation), 
we obtain Dhar' s formula \cite{dharo}:
\be
\E_{\mu} (N(x, y,\eta))=  G_{\mathbb{T}^q}(x,y),
\ee
where $G_{\mathbb{T}^q}$ denotes the Green's function of the \textit{Dirichlet random walk}
on $\mathbb{T}^q$. 

 The latter is defined as follows: at every vertex $x$ in $\mathbb{T}^q$, 
except the root, the random walk is killed with probability 
$(q+1-\deg(x))/(q+1)$ and is moving to every neighbouring vertex 
with probability $1/(q+1)$. At the root, the Dirichlet 
random walk is killed with probability
$(q-\deg(o))/q$ and is moving to a descendant with probability 
$1/q$.

By the identification between branching processes and $q$-ary trees with traps
 described at the beginning of this section, the Green's function  $G_{\mathbb{T}^q}$  of 
the Dirichlet random walk on $\mathbb{T}^q$ is  in distribution  equal to the Green's function of 
the random walk on $\mathcal{T}_q$ killed upon hitting a trap  (for a given trap configuration $\omega$,
with trapping 
probabilities $p_x=p$ for all $x\in \mathcal{T}_q$,  $x\not=o$, and the root $o$ is trapless), 
i.e. the Green's function 
 $\mathcal{G}_{T(\omega)}$ which is defined in  \eqref{def:eq_G} and estimated in 
 Proposition \ref{diam}.  We denote by $\mathbb{P}_{\mathbb{T}^q}$ 
the distribution over the realizations of the binomial tree
 and by $\mathbb{E}_{\mathbb{T}^q}$ the corresponding expectation.

This leads to the following  Theorem. \eec
\bt\label{Th:AvFinite}
We have for all $x\in \mathbb{T}^q$,

\begin{itemize}
\item[a)] annealed exponential bound on the diameter  of the avalanche: 
there is a constant $c:=c(q)>0$ such  that
\be 
\mathbb{P}_{\mathbb{T}^q}(\mu( ( \text{diam} (Av(x,\eta))>n))\leq e^{-cn}
\ee
\item[b)] and annealed finite expected avalanche size:
\be 
\mathbb{E}_{\mathbb{T}^q}(\mathbb{E}_{\mu}(|Av(x,\eta)|))<\infty.
\ee 
\
\end{itemize}

\et
\bpr
In part a),  notice that Markov's inequality and Dhar's formula give
\[
\begin{split}
\mu (y\in Av(x,\eta)) & =\mu (N(x,y ,\eta)\geq 1) \\
& \leq \mathbb{E}_{\mu}(N(x,y ,\eta)) = G_{\mathbb{T}^q} (x,y).
\end{split}
\]

Hence
\[
\begin{split}
\mathbb{P}_{\mathbb{T}^q}\left(\mu(diam(Av(x,\eta)) >n)\right)
 & =  \mathbb{P}_{\mathbb{T}^q}\left(\mu\left ( \exists y : d(x,y)>n, Av(x,\eta) \ni y\right)\right) \\
& \leq \mathbb{P}_{\mathbb{T}^q} \left ( \sum_{y: d(x,y)>n, \, y\in \mathbb{T}^q} \mu(Av(x,\eta) \ni y )\right ) \\
& \leq \mathbb{E}_{\mathbb{T}^q} \left ( \sum_{y: d(x,y)>n, \, y\in \mathbb{T}^q} G_{\mathbb{T}^q} (x,y)\right ) \\
\end{split}
\]
and the result follows by \eqref{diamgreen}.

For part b),  averaging over all tree realizations $\mathbb{T}^q$, using the equality in distribution of Green's functions and \eqref{sumgreen}, we get 
\[
\begin{split}
\mathbb{E}_{\mathbb{T}^q}(\mathbb{E}_{\mu}(|Av(x,\eta)|)) & = \mathbb{E}_{\mathbb{T}^q} 
\left(\sum_{y\in \mathbb{T}^q} \mu (y\in Av(x,\eta)) \right) \\
& \leq  \mathbb{E}_{\mathbb{T}^q} \left(\sum_{y\in \mathbb{T}^q} 
G_{\mathbb{T}^q}(x,y) \right)\\
& =  E \left (\sum_{x \in \caT_q} \mathcal{G}_{T(\omega)} (x,y) \right )\\
& =  \mathbf{E}_x(T(\omega) ),
\end{split}
\]
which is finite as given in Proposition \ref{expobound}.
\epr
\subsection{Other random trees}

In this section we want to stress that Theorem \ref{Th:AvFinite} 
applies also to other related models, namely the sandpile model 
on the random binary tree as in \cite{rrs} and a general Galton-Watson 
tree with bounded offspring distribution.\\

\textbf{(i)} In the {\it random binary tree}, each vertex has either 
2 descendants with probability $p$  or no descendant  with probability $1-p$. 
The Green's function of the Dirichlet random walk on this random binary tree can be 
then dominated by the Green's function of a random walk on a full binary tree 
with traps defined as follows. Each vertex of the full binary tree   $\mathcal{T}_2$   is a trap 
with probability $p$, independently of  other vertices.  The trap is effective with 
probability ${2}/{3}$, i.e., the random walk is killed with probability ${2}/{3}$
 upon hitting a trap. Then the Green's function of the latter random walk is 
 dominating  the Green's function of the Dirichlet random walk
  since the traps are not perfect and the random walk can survive 
 upon hitting a trap.\\

\textbf{(ii)}  A \textit{non-homogeneous branching process} 
$\mathbb{T}^q$ starts from a single individual. It has offspring probabilities
$(p_x(k))_{k\in \{0,..,q\}, x \in \mathbb{T}^q}$ to have $0,...,q$ descendants 
(so they possibly depend on the vertex $x$),  
such that $p_x(k) >p >0$ uniformly in $k$ and $x$. The Dirichlet's random walk 
Green's function can be again  dominated  by a trapped random walk on the 
full $q$-ary tree $\mathcal{T}_q$ as follows. Every vertex $x\in \mathcal{T}_q$ 
is a trap with probability $1-p_x(q)$. Upon hitting a trap, the random walk is 
killed with probability ${1}/{(q+1)}$.

\section{The  Abelian sandpile model with additional  dissipative sites 
 on the lattice $\Zd$}\label{sec:diss}
In this section, we consider the sandpile model on the
lattice $\Zd$ and add  to it  dissipative sites. More precisely, for the finite box
\be\label{finitebox}
\Lambda_n: =[-n,n]^d\cap \Zd,
\ee
 we consider the toppling matrix indexed by sites $x,y\in \Lambda_n$ and additionally parametrized by
a set $D\subset \Zd$ (finite or infinite) of  dissipative sites. 
We further denote $D_n:=D\cap \Lambda_n$ and $D^c_n:= D^c \cap \Lambda_n$,
\be\label{distopma}
\Delta_{x,y}^{D_n}
=
\begin{cases}
-1\ \ & \text{for}\ x,y\in \Lambda_n, x\sim y\\
2d+1\ &\text{for}\ x=y, \, x \in D_n\\
2d\  &\text{for}\ x=y, \, x \in D^c_n. 
\end{cases}
\ee
As before, we denote by $\caS_n$ stabilization with toppling matrix \eqref{distopma} within $\Lambda_n$
and $Av(x,\eta)$ the set of sites of $\Lambda_n$ which have to be toppled at least once upon 
addition of one grain at $x$. If $D=\emptyset$ then we will simply write 
$\Delta^{D_n}=\Delta^{\Lambda_n}$.
Furthermore, we call $\mu_n$ the uniform measure on recurrent configurations 
corresponding to the toppling matrix \eqref{distopma} and by $\mu$ as before the weak limit
\[
\mu = \lim_{n\rightarrow \infty}\mu_n,
\]
which exists, see, for example, \cite{jaraio}. 
\begin{definition}\label{Def:nonC}
We call the sandpile model  with dissipative sites $D$ on $\mathbb{Z}^d$ {\rm non-critical} 
if  for all $x\in \Z^d$ 
\be\label{nccond}
\limsup_{n\rightarrow \infty}\E_{\mu_n} (|Av(x,\eta)|)<\infty.
\ee
Otherwise, we call it {\rm critical}.
\end{definition}
Whether the model is critical or not will of course depend on the choice of the set of dissipative sites. 
Let us comment on Definition \ref{Def:nonC}. The sandpile model was introduced as a toy model  
displaying self-organized criticality, which is characterized by power-law behaviour 
of certain quantities like the avalanche distribution or two-point correlation functions. 
It is known in the mean-field setting \cite{MF}
and on homogeneous trees \cite{dhar} 
 that  we have the asymptotics 
\begin{equation}\label{muMF}
 \mu_n(|Av(0,\eta)|>k) \approx k^{-1/2}, \qquad \text{ as } k,n \rightarrow \infty.
\end{equation}
This  is also conjectured to hold above the critical dimension $d\geq 5$ \cite{priez}. 
A consequence of \eqref{muMF} is
\be\label{nccond}
\limsup_{n\rightarrow \infty}\E_{\mu_n} (|Av(x,\eta)|)=\infty.
\ee
Our definition is inspired by the analogous situation in percolation theory, 
where precisely at criticality the cluster is finite with probability 1 
and has infinite expectation,  while  in the sub-critical regime it has finite 
expectation; see chapter 1 in \cite{Gri}.
We define the finiteness of the avalanche cluster cardinality 
as the signature of non-criticality.
In what follows, we will be studying sufficient criteria on the set of 
dissipative sites ensuring that the model is no longer critical.\\

We will first characterize non-criticality by the Green's function  
associated to the toppling matrix and give sufficient conditions.  
Let us denote $G_n= \left(\Delta^{D_n}\right)^{-1}$. 
Note that for any given $x \in \Z^d$, there exists $n_0$ such that for all $n\geq n_0$, $x\in \Lambda_n$,
so that $ G_n(x,y)$ makes sense for $y \in \Z^d$. 
\bt\label{Thm:nonC}
Consider a sandpile model on $\Lambda_n$ with a set $D\subset{\mathbb{Z}^d}$ 
of dissipative sites. Then the model is non-critical if either
condition a) or b) is satisfied, and critical if condition c) is satisfied:
\begin{itemize}
\item[a)]   For all $x \in \Z^d$  
\be\label{sufco}
\limsup_{n\rightarrow \infty}\sum_{y\in \Lambda_n} G_n(x,y) <\infty.
\ee
\item[b)] For all $x \in \Z^d$, 
 the infinite volume Green's function
$
G(x,y)= \lim_{n\to\infty} G_n(x,y) 
$
exists.  Furthermore we have for all $x\in \mathbb{Z}^d$,
\[
\sum_{y\in \mathbb{Z}^d} G(x,y)<\infty.
\]
\item[c)] 
For all $x \in \Z^d$,  
$\lim_{n\to\infty} G_n(x,y)= G(x,y)$ is well-defined and there exists a  dissipative site $z$ such that
\[
\sum_{y\in \mathbb{Z}^d} G(z,y)=\infty.
\] 
\end{itemize}
\et
\bpr
This follows directly from the upper and lower bounds of Theorem 6.1 b) in \cite{JRS}.
\epr

Recall that as before we will use that there is a particular random walk associated 
to the Green's function, hence we can characterize non-criticality of the sandpile 
model via \bwc expected values of \ewc survival times of a corresponding random walk. Theorem \ref{sufthm} 
below describes sufficient
conditions for non-criticality of the dissipative sandpile model in terms of 
corresponding random walk estimates. 

Given a set of dissipative sites $D\subset \mathbb{Z}^d$ define a random walk 
$\widehat{X}=(\widehat{X}_k)_{k\in \mathbb{N}}$ on $\mathbb{Z}^d$ starting at  some point
 $x\in \mathbb{Z}^d$  as follows.  
If the random walk is at a non-dissipative site $x$, then it moves in the next step 
to one of its neighbours with probability ${1}/{(2d)}$. In the other case, it moves 
to one of its neighbours with probability ${1}/{(2d+1)}$ or is killed with 
probability ${1}/{(2d+1)}$.  
The  simple random walk on $\mathbb{Z}^d$ will be denoted by 
 $X=({X}_k)_{k\in \mathbb{N}}$ ,  $X_0=x$  with corresponding
 expectation $\E_x(\cdot)$. Call 
\[
\tau_x(D) =\inf \{ k\geq 1: X_k \in D \}
\]
the hitting or return time of the simple random walk to $D$.  
We denote by  $l_k(x)$ the corresponding local time,  i.e.,
the number of visits to $x$ of the simple random walk before time $k$. 
The following Theorem provides a number of sufficient conditions for non-criticality.
\begin{theorem}\label{sufthm}
 Let $\widehat{T}$ denote the survival time of the random walk 
 $(\widehat{X}_k)_{k\in \mathbb{N}}$ defined above with respect to some set 
 of dissipative sites $D\subset \mathbb{Z}^d$. Then the model is 
 non-critical if either of the following conditions is satisfied:
\begin{itemize}
\item[a)] For all $x\in \mathbb{Z}^d$, we have that $\mathbb{\widehat{E}}_x(\widehat{T})<\infty$.
\item[b)] For all $x\in \mathbb{Z}^d$, 
\[
\E_x \left ( \sum_{k=0}^{\infty} \prod_{z\in D} \left ( \frac{2d}{2d+1}\right )^{l_k(z)}\right ) <\infty.
\]
\item[c)] $|D^c| < \infty$. 
\item[d)] $D$ is such that there exists $\varphi: \mathbb{N}\rightarrow \R$ 
such that $\sum_{k=0}^{\infty} \left( \frac{2d}{2d+1}\right)^{\varphi(k)} <\infty$
and for all $x\in \mathbb{Z}^d$
\[
\sum_{k=0}^{\infty} \mathbb{P}_x \left(l_k(D^c) \geq k-\varphi(k) \right) <\infty.
\]
\item[e)]   For all $x\in \Z^d: \E_x \left(\sup_{y\in \mathbb{Z}^d} \tau_y(D)\right)<\infty$.
\item[f)] There exists a constant  $R >0$ such that 
\[
\sup_{x\in \Z^d} \inf_{y\in D} |x-y| =R +1< \infty,
\]
i.e., every point in $\Z^d$ is at most at distance $R+1$  from a point in $D$.
\end{itemize}
\end{theorem}
\bpr
\textbf{ a)} 
For this proof, we  introduce $(\overline{X}^{\Lambda_n}_t)_{t\ge 0}$,
the continuous-time random walk  in $\Lambda_n$  jumping at rate $2d$,
starting at $x\in\Lambda_n$ and killed upon leaving $\Lambda_n$;
we denote by $\overline{\E}^{\Lambda_n}_x(\cdot)$  the corresponding expectation;
moreover we denote by $\overline{X}=(\overline{X}_t)_{t\ge 0}$ 
the 
continuous-time random walk on $\mathbb{Z}^d$ jumping 
at rate $2d$, by
 $\overline{\E}(\cdot)$  the corresponding expectation  and by $\overline{l}_t(z)=\int_0^t \1_{\{\overline{X}_s=z\}} ds$ 
the corresponding  local time in $z$. 

Let us fix the set of dissipative sites $D$ and consider the associated 
toppling matrix $\Delta^{D_n}$ defined as in \eqref{distopma}.  Note that
\begin{equation}\label{potential}
\begin{split}
 G_n & = (\Delta^{D_n})^{-1} \\
& = (\Delta^{\Lambda_n} + \1_{D_n} Id)^{-1} =: (\Delta^{\Lambda_n} + V_{D_n} Id)^{-1}, 
\end{split}
\end{equation}
where $Id$ denotes the identity matrix in $\Lambda_n$ and $V_{D_n}=\1_{D_n}$ can be interpreted as a {\it potential}.  Finally fix $x\in \mathbb{Z}^d$  and $\Lambda_n$ for $n$ 
large enough such that $x\in \Lambda_n$.
In order to prove non-criticality: by Theorem \ref{Thm:nonC}(a),  we have to show that
\[
\limsup_{n\rightarrow \infty} \sum_{y\in \Lambda_n} G_n(x,y) <\infty.
\] 
By the Feynman-Kac formula  and using the exponential distribution of the jump times of the continuous time random walk  we can write
\beq\label{bolastro1}
G_n(x,y)
&=& \int_0^\infty \left(e^{-t\Delta^{D_n}}\right)_{x,y} dt\nonumber\\
&=&
\int_0^\infty {\overline\E}^{\Lambda_n}_x\left( e^{-\int_0^{t} V_{D_n}({\overline X}^{\Lambda_n}_s) ds} 
\1_{\{{\overline X}^{\Lambda_n}_t=y\}}\right) dt \nonumber\\
&\leq &
\int_0^\infty {\overline\E}_x\left( e^{-\int_0^{t} V_D(\overline{X}_s) ds} 
\1_{\{\overline{X}_t=y\}}\right) dt \nonumber\\
& = & \int_0^\infty \overline{\E}_x\left( e^{- \sum_{z\in D} \overline{l}_t(z)} 
\1_{\{\overline{X}_t=y\}}\right) dt \nonumber\\
&=& \frac1{2d}\E_x \left(\sum_{k=0}^\infty\prod_{z\in D}
\left(\frac{2d}{2d+1}\right)^{l_k(z)} \1_{\{X_k=y\}}\right).
\eeq
Then we have
\be\label{Thenwehave}
 \E_x \left [ \prod_{z\in D}\left(\frac{2d}{2d+1}\right)^{l_k(z)} 
 \1_{\{ X_k=y\}}\right ]= \widehat{\E}_x \left(\1_{\{\widehat{X}_k=y\}} 
 \1_{\{\widehat{T}>k\}}\right). 
\ee
Summing over all $y \in \mathbb{Z}^d$ gives
\[
\sum_{y \in \mathbb{Z}^d} \sum_{k \in \mathbb{N}} \widehat{\E}_x \left( \1_{\{\widehat{X}_k=y\}} \1_{\{\widehat{T}>k \}} \right)= \sum_{k=0}^\infty\widehat{\E}_x \left( \1_{\{\widehat{T}>k\}} \right)\leq \widehat{\E}_x (\widehat{T}).
\]
Therefore,
\begin{eqnarray}
\limsup_n\sum_{y\in \mathbb{Z}^d} G_n(x,y)  \leq 
& \frac{1}{2d}\E_x \left [\sum_{k=0}^{\infty} \prod_{z\in D }\left(\frac{2d}{2d+1}\right)^{l_k(z)} \right ]
\label{eqA} \leq  \frac{1}{2d}\widehat{\E}_x (\widehat{T}).
\end{eqnarray}
The claim follows using Theorem \ref{Thm:nonC} a). \\
\textbf{b)} Condition b) implies that the r.h.s. of \eqref{eqA} 
is finite and hence again  the condition in Theorem  \ref{Thm:nonC} a) is satisfied.\\
\textbf{c)} We will show that condition c) implies condition b).  
It is enough to consider $x=o$. Call $l_k(G)=\sum_{x\in G} l_k(x)$ 
for $G\subset \mathbb{Z}^d$. For some $\alpha\in (0,1)$, we write
\beq
&&\E_o \left(\left(\frac{2d}{2d+1}\right)^{l_k(D)}\right)
\nonumber\\
&=&
\E_o \left(\left(\frac{2d}{2d+1}\right)^{l_k(D)} \1_{\{ l_k(D)>\alpha k \}}\right)
+
\E_o \left(\left(\frac{2d}{2d+1}\right)^{l_k(D)} \1_{ \{l_k(D)\leq \alpha k \} }\right)
\nonumber\\
&\leq &
\left(\frac{2d}{2d+1}\right)^{\alpha k}
+
\E_o \left(\1_{\{ l_k(D) \leq \alpha k\}}\right).
\eeq
It now suffices to estimate the probability
$\pee_o \left( l_k(D) \leq \alpha k\right)$
and show that it is summable in $k$.
Notice that because $l_k(D)+ l_k(D^c)=k$ this amounts to estimate the probability
\[
\E_o \left(\1_{ \{l_k(D^c) \geq (1-\alpha) k)\} }\right).
\]
When $D^c$ is finite, this is simple.
We have the following bound on local time tails from 
Lemma 1 section 3 of \cite{holla}: there exist $a, b>0$ such that
for all $\delta>0$
\[
\pee_o \left(\sup_{x\in \mathbb{Z}^d} l_k(x) > k^{1/2 + \delta} \right)\leq a e^{-b k^{\delta/2}}.
\]
As a consequence, using
$l_k(D^c)< |D^c|\sup_{x\in \mathbb{Z}^d} l_k(x)$ we obtain for $\delta=\frac{1}{2}$,
\be\label{biri}
\pee_o\left(l_k(D^c) \geq (1-\alpha)k \right)\leq \pee\left(\sup_{x\in \mathbb{Z}^d} l_k(x) \geq \frac{(1-\alpha)k}{|D^c|}\right)\leq a'e^{-b' k^{1/4}},
\ee
which is summable in $k$. Non-criticality follows now from arguments in part a).\\
\noindent
\textbf{d)} This follows from the proof in c) by replacing  $\alpha k$ by $\varphi(k)$. \\
\noindent
\textbf{e)} Let us call $T_1, T_2,...$ the successive hitting times of the set $D$ of the simple random walk $X$. Every time that $D$ is visited, the corresponding trapped walk $\widehat{X}$ is killed with probability $\frac{1}{2d+1}$. The survival time $\widehat{T}$ of the killed walk $\widehat{X}$ starting from $x$ is a sum,
\[
\widehat{T}  \leq \sum_{i=1}^N \tau_{X_{T_i}} (D),
\]
where $N$ is a geometric random variable with parameter $\frac{1}{2d+1}$, 
independent of the simple random walk $X$. 
Therefore,

\[
\widehat{\E}_x(\widehat{T})\leq (2d+1) \E_x \left(\sup_{y\in \mathbb{Z}^d} \tau_y(D)\right) <\infty,
\]

then the claim follows from a).\\
\textbf{f)}
Denote for $x\in\Zd$ the ball with radius $R$,
\[
B(x,R):=\{ y\in \Zd: |x-y|\leq R\}.
\]
Upon exiting $B(x,R)$,
there is a strictly positive probability that a point of $D$ is hit. This probability is bounded from below
by a number $\kappa_R$ only depending on $R$. Indeed denote $\si_{B(x,R)}$  the exit time of $B(x,R)$.   Note that by translation invariance the distribution of $\sigma$ does not depend on $x$. Then
we can choose $\kappa_R= \inf_{z\in B(x,R)}\inf_{y} \pee_z( X_{\si_{B(x,R)}}=y)$. If no point of $D$ is
hit upon existing $B(x,R)$, then start from the exit point and look at the exit time of the ball with radius $R$
from that point.  We now see that from $x$
the hitting time of $D$ is bounded from above by a geometric sum of
exit times of balls of radius $R$ of which the distribution does not depend on $x$.
Therefore the condition of item e) is satisfied.
\epr

\br
One can also perform the sum over $y$ in the continuous-time expression appearing in the third line of   \eqref{bolastro1}.
This leads to the sufficient criterion of non-criticality
\begin{equation}\label{PAM}
\int_0^\infty \overline{\E}_x\left( e^{-\int_0^{t} V_D(\overline{X}_s) ds} \right)\ dt<\infty
\end{equation}
 for all $x\in \Z^d$. 
Notice that this corresponds to the total mass in the parabolic Anderson model
\cite{konig}, integrated over time. The parabolic Anderson model is a random walk model
in a random potential $V$. 
It also corresponds to the expected survival
time of a continuous-time random walk which is trapped at rate $1$ on dissipative sites and
moves to every neighboring site at rate one.
\er
The sufficient conditions for non-criticality provided 
in items c) - f) of Theorem \ref{sufthm} are not necessary. 
We provide two further examples illustrating critical versus non-critical
 behavior as a function of the ``size'' of the set of dissipative sites which are
not covered  by items c) - f) of Theorem \ref{sufthm}.

\textbf{(i)} First, if we make all sites of the $x$-axis of the 
two-dimensional square lattice dissipative, and all other
sites ordinary, then the model is critical because the expected hitting 
time of the $x$-axis of a two-dimensional simple random walk started
from a point  outside  the $x$-axis is infinite. By the same argument, 
if the set of dissipative sites is a lower dimensional subset
of $\Zd$ (such as a hyperplane intersected with $\Zd$), the model is still critical.

\textbf{(ii)} As a second example, consider the sandpile model on $\Z^2$ where 
we put dissipative sites on
a sequence of horizontal lines $y=0, y=-r_1, y=r_1, y=-r_2, y=r_2,\ldots$,
 where the distances $r_n$ form an increasing sequence
such that the gaps $r_{n+1}-r_n$ diverge in the limit $n\to\infty$.
Then the expected survival time starting at the origin is bounded by
\[
\frac14\sum_{i=1}^N (r_{i+1}-r_i)^2,
\]
where $N$ is an independent geometric random variable with success probability 
$1/(d+1)$ (see also the proof of Theorem \ref{sufthm} e). This can be seen from the fact that
the expected  hitting time of the set  $\{a,b\}$  of simple random walk starting 
at $x\in (a,b)$ is bounded by $(b-a)^2/4$.
Therefore, as long as the sum
\[
\sum_{k=0}^\infty \left(\frac{d}{d+1} \right)^k  (r_{k+1}-r_k)^2
\]
is convergent, the model is non-critical. 
Indeed one can choose the distances between two successive lines $r_{k+1}-r_k$  to grow faster than any polynomial in $k$. So this includes cases where the set of dissipative sites does not have a positive density in $\Z^2$ and hence cases not 
covered by the sufficient conditions provided in items c) - f) of Theorem \ref{sufthm}. 

\section{The Abelian avalanche model with sinks and sources}\label{sec:disssource}
In this section, we study a model with sources and sinks. Sources will be sites 
where upon toppling mass is gained, whereas sinks are sites where upon toppling 
mass is lost. To study such a system,  it is more convenient
 to  work in  the continuous-height setting  and therefore put ourselves in the context 
 of the Abelian avalanche model  (which was introduced in \cite{G93},
 then studied in \cite{JRS}),  which is the natural continuous-height 
 counterpart of the Abelian sandpile model. 

More precisely, we will consider a model which has dissipative sites in a 
set $D\subset\Zd$ and source sites in a set $S\subset\Zd$ and in which
the amount of mass transferred to neighbors  upon a toppling  is governed by a parameter $\gamma>0$, 
the amount of mass lost
upon toppling a (bulk) dissipative site is  governed by a parameter 
$\alpha$, and the amount of mass gained 
upon toppling a source site is  governed by a parameter  $\beta$.
This is what we call the \textit{Abelian avalanche model with sources and sinks}, 
with parameters $(D,S,\alpha, \beta,\gamma)$.  Note that one can retrieve the Abelian sandpile model from the Abelian avalanche model by $\gamma, \beta \rightarrow 0$ as in \cite{JRS}.

To define such a model, we first put ourselves in finite volume, and
 we  consider the following toppling matrix indexed by sites  $x,y\in\Lambda_n=[-n,n]^d\cap\Zd$,
 with $D_n=D\cap\Lambda_n, S_n=S\cap\Lambda_n$, 
\be\label{dissoutopma}
\Delta_{x,y}^{D_n, S_n, \alpha, \beta}
=
\begin{cases}
-\gamma\ \ &\text{for}\ x,y\in \Lambda_n, x\sim y\\
2d\gamma+\alpha\ &\text{for}\ x=y, \, x \in D_n\\
2d\gamma\ &\text{for}\ x=y, \, x\in \Lambda_n \setminus (D_n \cap S_n)\\
2d\gamma-\beta\ &\text{for}\ x=y, \, x\in S_n. 
\end{cases}
\ee
This defines a continuous-height sandpile model, called the 
{\it Abelian avalanche model} with
dissipative sites in $D_n$ and source sites in $S_n$.
The interpretation of the toppling matrix is, as already announced before, 
that a site is stable  if its height is  respectively below $2d\gamma$ for an 
``ordinary'' site $x\in \Lambda_n\setminus(D_n\cup S_n)$, below
$2d\gamma +\alpha$ for a dissipative site $x\in D_n$, and below 
$2d-\beta$ for a source site $x\in S_n$.

Upon toppling of a site $x\in\Lambda_n$, a mass $\gamma$ is transferred 
to each neighbor  of $x$  in $\Lambda_n$. 
This means that upon toppling, for a dissipative site which is not on the boundary, 
mass $\alpha$ is  {\em lost}, whereas for a source site
which is not on the boundary,
mass $\beta$ is {\em gained}.

For a height configuration $\eta: \Lambda_n \rightarrow [0,\infty)$ we define 
as before $\mathcal{S}_n(\eta)$ to be its stabilization according to the toppling 
matrix \eqref{dissoutopma}, provided $\eta$ is stabilizable, i.e. provided 
there exists a sequence of legal topplings with stable final result. 
If this is the case, then by the same argument as in  \cite{redigo}, 
stabilization is unique and well-defined.
We therefore assume first that $D, S, \alpha, \beta,\gamma $ are chosen in such a way 
that all $\eta$ are stabilizable for the toppling matrix
\eqref{dissoutopma} for all $n\in\N$. We then say that the model with parameters 
$(D, S, \alpha, \beta,\gamma)$ is {\em well-defined}.
 
In that case, we are in the same setting as the standard Abelian avalanche model, 
i.e., there exists a unique stationary measure
$\mu_{n}$ which is the uniform measure on recurrent configurations (i.e., 
configurations which are ``burnable'' -- through the burning algorithm, see
\cite{JRS}) and we have Dhar's formula
\[
\mathbb{E}_{n}(N(x,y, \eta)) = G_{n}(x,y),
\]
where $G_{n}(x,y)$ is again $(\Delta^{D_n,S_n,\alpha,\beta})^{-1}_{x,y}$ 
and where $N(x,y, \eta)$ denotes the number of topplings at $y$ needed
to stabilize $\eta+\delta_x$. As before, we then define the infinite volume model 
to be \textit{non-critical}  as in Definition \ref{Def:nonC}.
A sufficient condition for non-criticality is then (cf. \eqref{sufco})
\[
\limsup_{n\rightarrow \infty} \sum_{y\in \Lambda_n}G_{n}(x,y) <\infty.
\]

In analogy with the ``potential'' $V_D$ defined in \eqref{potential}, in our setting, we define the potential
\begin{equation}\label{sourcepot}
V_{D,S}(x)=
\begin{cases}
+\alpha\ &\text{if}\ x  \in D\\
0\ &\text{if}\ x\ \text{ordinary}\\
-\beta\ &\text{if}\ x \in S.
\end{cases}
\end{equation}
In analogy with \eqref{PAM}, we then have the following sufficient criterion 
of non-criticality of the infinite-volume model with parameters $(D, S, \alpha, \beta, \gamma)$.
\begin{proposition}
Provided the model with parameters  $(D,S,\alpha, \beta,\gamma)$ is well-defined, 
it is non-critical if for all $x\in\Zd$ we have
\be\label{genoncrit}
\int_0^\infty \overline{\E}_{\gamma, x}\left( e^{-\int_0^{t} V_{D,S}(\overline{X}_s) ds} \right) dt<\infty,
\ee
where $\overline{\E}_{\gamma, x}(\cdot)$ is the expectation w.r.t. the 
continuous-time random walk  $(\overline{X}_s)_{s\ge 0}$ with rate $2d\gamma$ starting at $x$.
\end{proposition}
The expectation $\overline{\E}_{\gamma, x}\left( e^{-\int_0^{t} V_{D,S}(\overline{X}_s) ds} \right)$ 
can be interpreted as the total mass at time $t>0$ starting
from a unit mass at time zero which is splitting at rate $\beta$ (in two 
unit masses) on source sites and killed at rate $\alpha$ at
dissipative sites, and besides is moving according to continuous-time random 
walk at rate  $2d\gamma$. Notice that
\begin{equation}\label{CRW}
\overline{\E}_{\gamma, x}\left( e^{-\int_0^{t} V_{D,S}(\overline{X}_s) ds}\right)
= \overline{\E}_x \left( e^{-\frac1\gamma\int_0^{t\gamma} V_{D,S}(\overline{X}_s)ds}\right),
\end{equation}
where we remind the reader that $\overline{\mathbb{E}}$ denoted the expectation w.r.t. the continuous time random walk. 
We will look at a finite number of source sites  and everywhere else dissipative sites, 
and show that for $\gamma$ large enough, the
model is \textit{not critical}. 
\begin{theorem}
Let $S$ be finite and $D:=\mathbb{Z}^d\setminus S$. Then
for $\gamma$ large enough, the sandpile model with toppling matrix
\eqref{dissoutopma} is not critical.
\end{theorem}
\bpr
We start with a single source site at the origin, $S=\{o\}$.
In that case
we have using \eqref{CRW}
\be\label{bimbi}
\overline{\E}_{\gamma, o}\left( e^{-\int_0^{t} V_{D,S}(\overline{X}_s)ds} \right)= \overline{\E}_o\left( e^{-\alpha t} e^{\frac{\alpha+\beta}{\gamma}l_{t\gamma}(\{o\})}\right)
\ee
where as before $l_t(G)=\int_0^t \1_{\{\overline{X}_s\in G\}} \, ds$ denotes the local time
associated to the random walk $(\overline{X}_s)_{s\ge 0}$. 
Denote
\be\label{free}
F(m)= \lim_{t\to\infty}\frac1t\log \overline{\E}_o \left ( e^{m\, l_{t}(\{o\}) } \right ).
\ee
Notice that this is exactly the free energy of the homogeneous pinning model; see \cite{giacomin}. The homogeneous pinning model is a  random polymer model which is penalized or rewarded upon touching the site $o$. 
 From  \cite{giacomin} Theorem 2.10,  we conclude that around $m\approx 0$ the behavior is
\[
F(m) = O(m^2)
\]
in $d=1, 3$.
In particular,
\be\label{bolastro}
\lim_{m\to 0} \frac1m F(m)=0.
\ee
Therefore,
\[
\overline{\E}_o \left ( e^{\frac{\alpha+\beta}{\gamma}l_{t\gamma}(\{o\})}\right )
\approx e^{\left(t\gamma F\left(\frac{\alpha+\beta}{\gamma}\right)+ o(t)\right)}
\]
by \eqref{bolastro},
and as $\gamma\to\infty$, we have
\be\label{freelim}
\gamma F\left(\frac{\alpha+\beta}{\gamma}\right)\to 0\ \text{as}\ \gamma\to\infty.
\ee
As a consequence, the right hand side of \eqref{bimbi} is integrable as a function of $t$ for $\gamma$ 
large enough.
In $d=2$, \eqref{bolastro} still holds. In $d\geq 4$, $F(m)=0$ for $m\in [0,m_c)$ 
with $m_c>0$; so also in that case
the rhs of \eqref{bimbi} is integrable for $\gamma$ large enough.

Finally, if we have a finite number of source sites $S$, then we need to estimate
\[
\overline{\E}_o \left ( e^{\frac{\alpha+\beta}{\gamma}l_{t\gamma}(S)}\right ),
\]
which by iteratively using Cauchy-Schwarz inequality can be estimated by
\[
e^{\gamma t F\left(2^{n-1} \frac{\alpha+\beta}{\gamma} \right) +o(t)},
\]
where $n=|S|$ is the number of source sites. The result then follows again from 
\eqref{freelim} as before.
\epr

{\bf Acknowledgements.} We thank 
MAP5 lab at Universit\'e Paris Descartes
and Delft University of Technology for financial support and hospitality.

\end{document}